\documentclass[11pt]{article}
\usepackage{mathrsfs}
\usepackage{amsfonts}
\usepackage{amssymb}
\usepackage{amsmath,amscd}
\usepackage[dvips]{graphicx}
\usepackage[all]{xy}
\usepackage{graphicx}
\usepackage{fancyhdr}
\usepackage{mathptm,pslatex}
\usepackage{amsthm}

\begin{document}

\sloppy
\newtheorem{Def}{Definition}[section]
\newtheorem{Bsp}{Example}[section]
\newtheorem{Prop}[Def]{Proposition}
\newtheorem{Theo}[Def]{Theorem}
\newtheorem{Lem}[Def]{Lemma}
\newtheorem{Koro}[Def]{Corollary}
\theoremstyle{definition}
\newtheorem{Rem}[Def]{Remark}

\newcommand{\add}{{\rm add\, }}
\newcommand{\fd}{{\rm fin.dim\, }}
\newcommand{\gd}{{\rm gl.dim\, }}
\newcommand{\dd}{{\rm dom.dim\, }}
\newcommand{\E}{{\rm E}}
\newcommand{\Mor}{{\rm Morph\, }}
\newcommand{\End}{{\rm End\, }}
\newcommand{\ind}{{\rm ind\,}}
\newcommand{\rsd}{{\rm res.dim\,}}
\newcommand{\rd} {{\rm rep.dim\, }}
\newcommand{\ol}{\overline}
\newcommand{\overpr}{$\hfill\square$}
\newcommand{\rad}{{\rm rad\,}}
\newcommand{\soc}{{\rm soc\,}}
\renewcommand{\top}{{\rm top\,}}
\newcommand{\pd}{{\rm pd\, }}
\newcommand{\id}{{\rm inj.dim\, }}
\newcommand{\Fac}{{\rm Fac\,}}
\newcommand{\DTr}{{\rm DTr\, }}
\newcommand{\cpx}[1]{#1^{\bullet}}
\newcommand{\D}[1]{{\mathscr D}(#1)}
\newcommand{\Dz}[1]{{\mathscr D}^+(#1)}
\newcommand{\Df}[1]{{\mathscr D}^-(#1)}
\newcommand{\Db}[1]{{\mathscr D}^b(#1)}
\newcommand{\C}[1]{{\mathscr C}(#1)}
\newcommand{\Cz}[1]{{\mathscr C}^+(#1)}
\newcommand{\Cf}[1]{{\mathscr C}^-(#1)}
\newcommand{\Cb}[1]{{\mathscr C}^b(#1)}
\newcommand{\K}[1]{{\mathscr K}(#1)}
\newcommand{\Kz}[1]{{\mathscr K}^+(#1)}
\newcommand{\Kf}[1]{{\mathscr  K}^-(#1)}
\newcommand{\Kb}[1]{{\mathscr K}^b(#1)}
\newcommand{\modcat}{\ensuremath{\mbox{{\rm -mod}}}}
\newcommand{\Modcat}{\ensuremath{\mbox{{\rm -Mod}}}}

\newcommand{\stmodcat}[1]{#1\mbox{{\rm -{\underline{mod}}}}}
\newcommand{\pmodcat}[1]{#1\mbox{{\rm -proj}}}
\newcommand{\imodcat}[1]{#1\mbox{{\rm -inj}}}
\newcommand{\Pmodcat}[1]{#1\mbox{{\rm -Proj}}}
\newcommand{\op}{^{\rm op}}
\newcommand{\otimesL}{\otimes^{\rm\bf L}}
\newcommand{\rHom}{{\rm\bf R}{\rm Hom}\,}
\newcommand{\Hom}{{\rm Hom}}
\newcommand{\coker}{{\rm Coker}}
\renewcommand{\ker}{{\rm Ker}}
\newcommand{\Img}{{\rm Im}}
\newcommand{\ext}{{\rm Ext}}
\newcommand{\StHom}{{\rm \underline{Hom} \, }}
\renewcommand{\rad}{{\rm rad}}

\def\vez{\varepsilon}
\def\bz{\bigoplus}
\def\sz {\oplus}
\def\epa{\xrightarrow}
\def\inja{\hookrightarrow}

\newcommand{\lra}{\longrightarrow}
\newcommand{\lraf}[1]{\stackrel{#1}{\lra}}
\newcommand{\ra}{\rightarrow}
\newcommand{\dk}{{\rm dim_{_{k}}}}
\newcommand{\colim}{{\rm colim\, }}
\newcommand{\limt}{{\rm lim\, }}
\newcommand{\Add}{{\rm Add\, }}
\newcommand{\Tor}{{\rm Tor}}
\newcommand{\Cogen}{{\rm Cogen}}

{\Large \bf
\begin{center}
Finitistic dimension conjecture and radical-power extensions\\
{{\small \it Dedicated to Claus Micheal Ringel on the occasion of his 70$^{th}$ birthday}}
\end{center}}

\medskip
\medskip
\centerline{\textbf{Chengxi Wang} and \textbf{Changchang Xi}$^*$}

\renewcommand{\thefootnote}{\alph{footnote}}
\setcounter{footnote}{-1} \footnote{ $^*$ Corresponding author.
Email: xicc@cnu.edu.cn; Fax: 0086 10 68903637.}
\renewcommand{\thefootnote}{\alph{footnote}}
\setcounter{footnote}{-1} 
\renewcommand{\thefootnote}{\alph{footnote}}
\setcounter{footnote}{-1} 

\begin{abstract}
The  finitistic dimension conjecture asserts that any finite-dimensional algebra over a field should have finite finitistic dimension. Recently, this conjecture is
reduced to studying finitistic dimensions for extensions of algebras. In this paper, we investigate those extensions of Artin algebras in which some radical-power of
smaller algebras is a one-sided ideal in bigger algebras. Our results, however, are formulated more generally for an arbitrary ideal: Let $B\subseteq A$
be an extension of Artin algebras
and $I$ an ideal of $B$ such that the full subcategory of $B/I$-modules is $B$-syzygy-finite. Then: (1) If the extension is right-bounded
(for example, $\pd(A_B)<\infty$),
$I\, A\,\rad(B)\subseteq B$ and $\fd(A)<\infty$, then $\fd(B)<\infty$. (2) If $I\, \rad(B)$ is a left ideal of $A$ and $A$ is torsionless-finite, then
$\fd(B)<\infty$. Particularly, if $I$ is specified to a power of the radical of $B$, then our results not only
generalize some ones in the literature (see Corollaries \ref{1.1} and \ref{1.2}), but also provide some completely new ways to detect algebras of finite finitistic dimensions.
\end{abstract}

\section{Introduction}
Let $A$ be an Artin algebra. The finitistic dimension of $A$ is defined to be the supremum of the projective dimensions of finitely generated
left $A$-modules having finite projective dimension. The finitistic dimension conjecture says that any Artin algebra should have finite finitistic
dimension. This conjecture was initially a question by Rosenberg and Zelinsky, published by Bass in a paper in 1960 (see \cite{Bass}), and has attracted many
mathematicians in the last $5.5$ decades. Among them is Maurice Auslander who ``is considered to be one of the founders of the modern aspects of the
representation theory of artin algebras" (see \cite[p.501]{reitensmalosolberg}). ``One of his main interests in the theory of artin algebras was the
finitistic dimension conjecture and related homological conjectures" (see \cite[p.815]{reitensmalosolberg}). The conjecture becomes
nowadays one of the main conjectures in the representation theory of algebras (see \cite[Conjecture (11), p.410]{ars}) and has not only close relations
with algebraic geometry and model structure (see \cite{Membrillo, IEA}), but also intimate connections with the solutions of several
other not-yet-solved conjectures such as strong Nakayama conjecture (see \cite{CF}), generalized Nakayama conjecture (see \cite{AR2}), Nakayama
conjecture (see \cite{Naka}), Wakamatsu tilting conjecture (see \cite[Chapter IV, p.71]{beligiannis-reiten}) and Gorenstein symmetry conjecture
(see \cite{Yam}). All of these conjectures would be valid if so would be the finitistic dimension conjecture. Several special cases for the conjecture
to be true are verified (see, for example, \cite{Auslander3, EHIS, Green, GKK, HLM, Igusa, Z, Zim, recIII}), but it is not yet fully resolved in general. Actually, up to the present time, not many practical methods, so far as we know, are available to detect algebras of finite finitistic dimensions. It seems necessary to develop some methods for testing finiteness of finitistic dimensions for general algebras or even for some concrete examples.

In the recent papers \cite{I, II}, the conjecture is reduced to comparing finitistic dimensions of a pair of algebras instead of focusing only on one single algebra. More precisely, the following two statements are proved to be equivalent for a field $k$:

(1) The finitistic dimension of any finite-dimensional $k$-algebra is finite.

(2) For any extension $B\subseteq A$ of finite-dimensional $k$-algebras such that $\rad(B)$, the Jacobson radical
of $B$, is a left ideal in $A$, if $A$ has finite finitistic dimension, then $B$ has
finite finitistic dimension.

Along this line, the conjecture is further reduced, by a different method, to extensions of algebras with relative global dimension $1$,
where the ground field is assumed to be perfect (see \cite{XiXu}). Thus it seems quite worthy to consider such kinds of extensions of algebras
and to bound the
finitistic dimensions of smaller algebras in terms of the ones of bigger algebras which we would like to take as simple as possible. Moreover,
this kind of considerations, philosophically, seems to make sense because algebras usually may have simpler homological or representation-theoretical
properties than their subalgebras do. For instance, every finite-dimensional algebra over a field can be regarded as a subalgebra of some full matrix
algebra, while the latter obviously has simple representation theory and homological properties. Also, a lot of other examples show that the philosophy
of controlling finitistic dimensions by extension algebras could work powerfully (see also the example at the end of the last section).

Now, let us just mention a couple of known considerations in this direction. In the sequel, we denote by $\gd(A)$ and $\fd(A)$ the global and finitistic dimensions of an algebra $A$, respectively; and by $\pd(A_{B})$ the projective dimension of the right $B$-module $A$.

(i) Let $B\subseteq A$ be an extension of Artin algebras such that rad$(B)$ is a left ideal in $A$. Then:

  $\qquad$ (a) If $\rad(A)=\rad(B)A$ and $\gd(A)\le 4$, then $\fd(B)<\infty$ (see \cite[Theorem 3.7]{II}).

  $\qquad$ (b) If $\pd(A_B)<\infty$ and $\fd(A)< \infty$, then $\fd(B)<\infty$ (see \cite[Corollary 1.4]{XiXu}).

\smallskip
(ii) Let $C\subseteq B\subseteq  A$ be a chain of Artin algebras such that $\rad(C)$ and $\rad(B)$ are left ideals in $B$ and $A$, respectively. If $A$ is representation-finite, then $\fd(C)<\infty$ (see \cite[Theorem 4.5]{I}).

\medskip
In this paper, we continue to explore and compare finitistic dimensions of extensions $B\subseteq A$ of Artin algebras.
Since it happens quite often that $\rad(B)$ itself may not be a left ideal in $A$ but some power of it is actually a left or right ideal in $A$, our main goal in this paper is to extend results for the case
that $\rad(B)$ is a left ideal in $A$ to a more general case that rad$^s(B)$ is a left (or right) ideal in $A$ for some positive integer $s$. Such kinds of extensions may
be called \emph{radical-power} extensions. Our results in this paper will generalize some ones on left-idealized extensions (see \cite{I, II, XiXu}), recover a result of Green-Zimmermann-Huisgen
(see Corollary \ref{koro6}) and provide somewhat
handy methods to detect algebras of finite finitistic dimensions, especially those obtained as subalgebras from representation-finite algebras.

We shall carry out our discussion for extensions $B\subseteq A$ in a broader context by studying an arbitrary ideal of $B$ rather than a power of the radical of $B$.
The technical problem we encounter, even in the case of a higher power of the radical, is that the higher syzygies of $B$-modules admit no longer $A$-module structures.
So a crucial ingredient for bounding projective dimensions used in \cite{I, II, XiXu} is missing. To circumvent this problem here, we first use certain submodules of
torsionless $B$-modules to get $A$-module structures, and then establish certain reasonable short exact sequences connecting $B$-syzygies
with the lifted $A$-modules. Finally, we utilize the Igusa-Todorov function in \cite{Igusa} to estimate upper bounds of projective dimensions.

To state our main results more precisely, let us first recall some definitions.

Let $A$ be an Artin algebra and $M_A$ be a right $A$-module. Following \cite{KK}, a non-negative integer $n$ is called a \emph{bound} on the vanishing of $\Tor_A(M_A,-)$ if, for any $A$-module $Y$, whenever $\Tor^A_p(M_A,Y)=0$ for $p$ sufficiently large then $\Tor^A_p(M_A,Y)=0$ for all $p\ge n+1$. Of course, if $n:=\pd(M_A)<\infty$, then $n$ is a bound on the vanishing of $\Tor^A(M,-)$. As in \cite{KK}, if the $n$-th syzygy of a projective resolution of $M_A$ is periodic, then $n$ is a bound on the vanishing of $\Tor^A(M,-)$. So a bound on the vanishing of $\Tor^A(M,-)$ generalizes the period and finiteness of a projective resolution of $M$.

An extension $B\subseteq A$ of Artin algebras is said to be \emph{right-finite} if $\pd(A_B)<\infty$. Similarly, we can define left-finite extensions. Now, we introduce a generalization of right-finite extensions. The extension
$B\subseteq A$ is called a \emph{right-bounded} extension if there is a bound on the vanishing of $\Tor^B(A_B,-)$.
So, right-finite extensions are right-bounded, but the converse is not true in general by the examples in \cite{KK}.

For $m\ge 0$ and a full subcategory $\mathcal{C}$ of $A\modcat$, the category of all finitely generated left $A$-modules, we denote by $\Omega^m({\mathcal C})$
the full subcategory of $A$-mod consisting of those $A$-modules that are either projective or direct summands of $m$-th syzygies of $A$-modules in $\mathcal C$. So $\Omega^0(A\modcat)=A\modcat$
and $\Omega^1(A\modcat)$ is just the category of torsionless $A$-modules. We say that $\mathcal C$ is $m$-\emph{syzygy-finite} if there are only finitely many non-isomorphic indecomposable modules
in $\Omega^m(\mathcal{C})$, and \emph{syzygy-finite} (or more precisely, $A$-\emph{syzygy-finite}) if there is some $m\ge 0$ such that $\mathcal C$ is $m$-syzygy-finite. Particulary, the algebra $A$ is \emph{representation-finite} (respectively, \emph{torsionless-finite}) if $A\modcat$ is $0$-syzygy-finite (respectively, $1$-syzygy-finite). For further information on torsionless-finite algebras, we refer the reader to \cite{ringel}.

Now, our main results can be stated as follows.

\begin{Theo}
Let $B\subseteq A$ be a right-bounded extension of Artin algebras and $I$ be an ideal in $B$ such that $I\, A \, \rad(B)\subseteq B$ and the full subcategory of $B/I$-modules is $B$-syzygy-finite
(for example, $B/I$ is representation-finite). If $\fd(A)<\infty$, then $\fd(B)<\infty.$ \label{1.1b}
\end{Theo}

As an immediate consequence of Theorem \ref{1.1b}, we get a result on radical-power extensions of Artin algebras.

\begin{Koro}
Let $B\subseteq A$ be a right-finite extension of Artin algebras. Suppose that there is an integer $s\ge 0$ such that $\rad^{s}(B)A\,\rad(B)\subseteq B$ and that $B/\rad^{s}(B)$ is representation-finite. If $\fd(A)<\infty$, then $\fd(B)<\infty$.
\label{1.1}
\end{Koro}

Here, we understand rad$^0(B)=B$. In case $s=0$ and $\rad(B)$ is a left ideal in $A$, Corollary \ref{1.1} coincides with \cite[Corolary 1.4]{XiXu} for finiteness of finitistic dimensions.
But in the case $s=1$, Corollary \ref{1.1} seems to be new, comparing it with the results in \cite{cowey, KK, wei, I, II, XiXu}.

As is known, $\fd(A)$ may differ from $\fd(A^{\op})$, where $A^{\op}$ stands for the opposite algebra of $A$. This means that the notion of finitistic dimensions does not have left-right
symmetry and suggests that a conclusion analogous to Corollary \ref{1.1} for left-finite extensions $B\subseteq A$ with $\rad^{s}(B)$ being a left (or right) ideal of $A$ might not exist.
However, in this case, we get similar results in which we do not assume that extensions are left-bounded.

\begin{Theo} Let $B\subseteq A$ be an extension of Artin algebras and $I$ be an ideal in $B$ such that $I \, \rad(B)$ is a left ideal in $A$ and the full subcategory of $B/I$-modules
is $B$-syzygy-finite. If $A$ is torsionless-finite (for example, representation-finite or hereditary), then $\fd(B)<\infty$.
\label{1.3}
\end{Theo}

If $I=B$ in Theorem \ref{1.3}, then we recover \cite[Theorem 3.1]{I}. Now, specifying $I$ to some power of the radical of $B$, we then get the following corollary.

\begin{Koro}
Suppose that $B\subseteq A$ is an extension of Artin algebras such that $\rad^{s}(B)$ is a left ideal of $A$ and that $B/\rad^{s-1}(B)$ is representation-finite for some integer $s\geq 1$. If $A$ is torsionless-finite,  then $\fd(B)<\infty$.
\label{1.2}
\end{Koro}

Note that Corollary \ref{1.2} covers \cite[Proposition 4.9]{I} if we take $s=1$. In case $s=2$, the algebra $B/\rad^{s-1}(B)$ is automatically representation-finite, and therefore
Corollary \ref{1.2} takes a simple form. In this case, Corollary \ref{1.2} seems to appear for the first time in the work. Since there is a plenty of extensions $B\subseteq A$ such
that $\rad(B)$ itself is neither a left nor a right ideal in $A$ but some of its powers is a left or right ideal in $A$, our corollaries are proper generalizations of
some results in \cite{I, II, XiXu}. Also, the arguments of proofs of our results here are different from the earlier ones, though the common idea
for all proofs is the use of the Igusa-Todorov function. Note also that, comparing with the results in \cite{wei}, we do not
impose any homological conditions (such as finiteness of projective dimension) on ideals or powers of the
radical of $B$ since such conditions seem to be strong for our homological questions.

The contents of this note are arranged as follows: In Section \ref{sect2}, we first fix some terminology and notation, and then recall some basic facts needed in later proofs.
In Section \ref{sect3}, we give proofs of all the above-mentioned theorems and
corollaries. As a byproduct of our results, we re-obtain a main result in \cite{Green} (see Corollary \ref{koro6}). This section ends with an example to explain how our results can be used. It seems that, for this example, no previously known results could be applied
to detect the finiteness of its finitistic dimension.

\section{Preliminaries\label{sect2}}

In this section, we shall briefly recall some basic definitions, fix notation,
and collect some known results that are needed in later proofs.

Let $A$ be an Artin $R$-algebra, that is, $R$ is a commutative Artin ring with identity and $A$ is a finitely generated unitary $R$-module. We denote by $\rad(A)$ the Jacobson radical of $A$. By a module we always mean a finitely generated left module. Given an $A$-module $M$, we denote by $\rad_A(M)$ and $\Omega_A(M)$ the Jacobson radical and the first syzygy of $M$, respectively, and by $\pd(_AM)$ the projective dimension of $M$. Recall that the global dimension of $A$, denoted by $\gd(A)$, is defined to be the supremum of projective dimensions of modules in $A\modcat$, and the finitistic dimension of $A$, denoted by $\fd(A)$, is the supremum of projective dimensions of all those $A$-modules $X$ in $A\modcat$ with $\pd(_AX)<\infty$. So $\fd(A)\le \gd(A).$

For a homomorphism $f: X\ra Y$ in $A\modcat$, we denote by $\ker(f)$ and $\Img(f)$ the kernel and image of $f$, respectively. If $g:Y\ra Z$ is another homomorphism in $A\modcat$, we write $fg: X\ra Z$ for the composite of $f$ with $g$. In this way, $\Hom_A(X,Y)$ becomes naturally a left $\End_A(X)$- and right $\End_A(Y)$-bimodule.

Given an additive full subcategory $\mathscr X$ of $A\modcat$ and a module $M$ in $\mathscr X$, we say that $M$ is an \emph{additive generator} for $\mathscr X$ if $\mathscr X$ = $\add(M)$, where $\add(M)$ is the additive full category of $A\modcat$ generated by $M$.

By an extension $B\subseteq A$ of algebras we mean a pair of Artin algebras $A$ and $B$ such that $B$ is a subalgebra of $A$ with the same identity.

To estimate projective dimensions of modules, Igusa-Todorov introduced a useful function $\Psi$ (see \cite{Igusa}) which is defined as follows: Let $K(A)$ be the quotient of the free abelian group generated by the isomorphism classes $[M]$ of modules $M$ in $A$-mod modulo the relations:

(i) $[Y]=[X]+[Z]$ if $Y\simeq X\oplus Z$, and

(ii) $[P]=0$ if $P$ is a projective $A$-module.

\noindent Thus $K(A)$ is a free abelian group with a basis consisting of isomorphism classes of non-projective indecomposable $A$-modules. By employing the noetherian property of the ring of integers, Igusa and Todorov define a function $\Psi$ on this abelian group, which takes values of non-negative integers and has the following fundamental properties.

\begin{Lem}\emph{\cite{Igusa} } \label{it} Let $X$ and $Y$ be $A$-modules. Then:

$(1)$ If $\pd(_AX)<\infty$, then $\Psi(X) =\pd(_AX)$. 

$(2)$ $\Psi(X)\leq\Psi(X\oplus Y)$ and $\Psi(\bigoplus^{n}_{j=1}X)=\Psi(X)$ for any natural number $n\ge 1$.

$(3)$ If $0\ra X\ra Y\ra Z\ra 0$ is an exact sequence in $A\modcat$ with $\pd(_AZ)<\infty$, then $\Psi(Z)\leq\Psi(X\oplus Y)+1$.
\end{Lem}

Now, we recall the following well-known homological facts.

\begin{Lem} \label{facts} $(1)$ If there is an exact sequence $0\ra X_s \ra \dots \ra X_1 \ra X_0 \ra X \ra 0$ in $A\modcat,$ then $\pd(_AX) \le s + \max\{\pd(_AX_i)\mid 0\le i\le s\}$.

$(2)$ Given an exact sequence $ 0 \ra X \ra Y \ra Z\ra 0$ of $A$-modules, there are induced another two exact sequences $$0 \lra \Omega_{A}(Z) \lra X\oplus P \lra Y \lra 0 \; \mbox{\;  and \; } \; 0\lra \Omega_A(Y)\lra \Omega_A(Z)\oplus P'\lra X\lra 0$$
with $P$ and $P'$ projective $A$-modules. They are the so-called syzygy shifting.
\end{Lem}

Next, we recall a result of Auslander (see \cite[Chapt. III, Sec. 2]{AusRepdim}), which describes the global dimensions of the endomorphism algebras of modules.

\begin{Lem}\label{repdim} Let $n\ge 2$ be an integer. Suppose that $A$ is an algebra and $M$ is an $A$-module. Then $\gd(\End(_AM))\le n$ if and only if,
for each $A$-module $X$, there is an exact sequence $$0\lra M_{n-2}\lra \cdots \lra M_1\lra M_0\lra X \lra 0$$ in $A\modcat$ such that all $M_j\in \add(M)$ and $\Hom_A(M,-)$ preserves the exactness of this sequence.
\end{Lem}

We also need the following result from \cite{II}.

\begin{Lem}
Let $B\subseteq A$ be an extension of algebras such that $\rad(B)$ is a left ideal of $A$. Then, for any $B$-module $X$ and integer $i\geq 2$, the module $\Omega^i_B(X)$ has an $A$-module structure and there is a projective $A$-module $P$ and an $A$-module $Y$ such that $\Omega^i_B(X) \simeq \Omega_A(Y)\oplus P$ as $A$-modules.
\label{Lem4}
\end{Lem}

\section{Proofs of main results\label{sect3}}

In this section, we give proofs of all results mentioned in the introduction and deduce some of their consequences.

\subsection{Proof of Theorem \ref{1.1b}}

To prove Theorem \ref{1.1b}, we first have to make some preparations. Let $M$ be a right $A$-module. For convenience, we say that
a left $A$-module $_AX$ satisfies the \emph{$\Tor^A(M,-)$-vanishing condition} if $\Tor_p^A(M,X)=0$ for $p$ sufficiently large. For instance, modules of finite projective dimension always satisfy this condition.
Moreover, if $\pd(M_A)<\infty$, then all modules satisfy the $\Tor^A(M,-)$-vanishing condition. Note that the next two lemmas are true for a
more general class of rings, though they are just stated for Artin algebras.

\begin{Lem} \label{3.1}
Let $A$ and $B$ be Artin algebras and $_AM_B$ be an $A$-$B$-bimodule. Suppose that $_AM$ is projective and $n$ is a bound on the vanishing of $\Tor^B(M_B,-)$. Then, for any $B$-module $X$ with the $\Tor^B(M_B,-)$-vanishing condition and for any  integer $m\ge 1$, we have an isomorphism of $A$-modules $$\Omega_{A}(M\otimes_{B}\Omega^{n+m }_{B}(X))\simeq\Omega^{m}_{A}(M\otimes_{B}\Omega^{n+1}_B(X)).$$
\end{Lem}

{\it Proof.}
Let $m\ge 1$ be an integer and $X$ be a $B$-module satisfying the $\Tor^B(M_B,-)$-vanishing condition. Then we consider the projective cover $f: P\ra \Omega_B^{n+m-1}(X)$ of $\Omega_B^{n+m-1}(X)$ and the canonical exact sequence
$$0 \lra \Omega_B^{n+m}(X)\lra P \lra \Omega_B^{n+m-1}(X) \lra 0.$$
Since $n$ is a bound of the vanishing of $\Tor^B(M_B,-)$ and $n+m>n$, we have $\Tor_1^{B}(M, \Omega_B^{n+m-1}(X))\simeq\Tor_{n+m}^B(M,X)=0$. Thus the following sequence of $A$-modules is exact:
$$\xymatrix@C=0.5cm{
0 \ar[r] & M\otimes_B\Omega_B^{n+m}(X) \ar[r]^{} & M\otimes_{B}P \ar[r]^{} & M\otimes_{B}\Omega_B^{n+m-1}(X) \ar[r] & 0 }.$$By assumption, the module $_AM$ is projective. This yields that the $A$-module $M\otimes_BP$ is also projective. Thus $M\otimes_{B}\Omega_{B}^{n+m}(X)\simeq \Omega_{A}(M\otimes_{B}\Omega_{B}^{n+m-1}(X))\oplus Q $ with $Q$ a projective $A$-module. It follows that
$\Omega_{A}(M\otimes_{B}\Omega_{B}^{n+m}(X))\simeq \Omega_{A}^{2}(M\otimes_{B}\Omega_{B}^{n+m-1}(X))$ for all $m\ge 1$. Consequently, we have $\Omega_{A}(M\otimes_{B}\Omega_{B}^{n+m}(X))\simeq \Omega_{A}^{m}(M\otimes_{B}\Omega_{B}^{n+1}(X))$ by repeating the foregoing isomorphism. $\square$

\medskip
The following lemma is a consequence of the derived functors of tensor functors.

\begin{Lem}\label{3.2}
Let $B\subseteq A$ be a right-bounded extension of Artin algebras with $n$ a bound on the vanishing of $\Tor^B(A_B,-)$. If a $B$-module $_BX$ satisfies the $\Tor(A_B,-)$-vanishing condition, that is, $\Tor^B_p(A_B,X)=0$ for sufficiently large $p$, then, for any positive integer $i\geq n+1$, the following sequence of $B$-modules is exact:
$$\xymatrix@C=0.5cm{0 \ar[r] & \Omega_{B}^{i}(X) \ar[r]^{} & A\otimes_{B}\Omega_{B}^{i}(X) \ar[r]^{} & (A/B)\otimes_{B}\Omega_{B}^{i}(X) \ar[r] & 0 }.$$
\end{Lem}

{\it Proof.} Tensoring the exact sequence $0\ra B\ra A\ra A/B\ra 0$
by $\Omega_B^i(X)$, we get the following exact sequence of $B$-modules:
$$0\lra\Tor_1^B(A,\Omega_{B}^{i}(X))\lra\Tor_1^B(A/B,\Omega_{B}^{i}(X))\lra \Omega_{B}^{i}(X)\lra A\otimes_B\Omega_{B}^{i}(X)\lra (A/B)\otimes_B\Omega_{B}^{i}(X)\lra 0. $$
Evidently, $\Tor_1^B(A/B,\Omega_{B}^{j}(X))\simeq \Tor_{j+1}^B(A/B,X)\simeq \Tor_{j+1}^B(A,X)$ for $j\ge 1.$
By assumption, $n$ is a bound of the vanishing of $\Tor^B(A_B,-)$ and $X$ satisfies the $\Tor^B(A_B,-)$-vanishing condition. Therefore we have $\Tor_{p}^B(A,X)=0$ for all $p\ge n+1$. Thus $\Tor_1^B(A/B,\Omega_{B}^{i}(X))\simeq \Tor_{i+1}^B(A,X)=0$ for $i\ge n+1$. Consequently, we get the desired exact sequence. $\square$

\medskip
{\bf Proof of Theorem \ref{1.1b}}:

Let $B\subseteq A$ be a right-bounded extension of Artin algebras and $I$ be an ideal in $B$ such that $I\, A \, \rad(B)\subseteq B$ and the full subcategory $(B/I)\modcat$ of $B\modcat$
is $s$-syzygy-finite (with respect to taking $B$-syzygies) for some integer $s\ge 0$. Further, we assume that $n$ is a bound on the vanishing of $\Tor(A_B,-)$, and define  $m:=\fd(A)<\infty$. Let $X$ be a $B$-module with $\pd(_BX)<\infty$. Then  $X$ satisfies the $\Tor^B(A_B,-)$-vanishing condition. Now,
we consider the $B$-module $Y:=\Omega_{B}^{m+n+1}(X)$, take a projective cover $\pi: P\ra Y$ of $Y$, and form the following exact commutative diagram of $B$-modules:
 $$\xymatrix{
 &0\ar[d]^{}&0\ar[d]^{}&0\ar[d]^{}& \\
0\ar[r]^{}& \Omega_{B}(Y)\ar[r]^{}\ar[d]^{}&  P\ar[r]^{\pi}\ar[d]^{}& Y\ar[r]^{}\ar[d]^{}& 0\\
0\ar[r]^{}&A\otimes_{B}\Omega_B(Y)\ar[r]^{}\ar[d]^{}&  A\otimes_{B}P \ar[r]^{}\ar[d]^{}& A\otimes_{B}Y\ar[r]^{}\ar[d]^{}& 0\\
\cdots\ar[r]^{}& (A/B)\otimes_{B}\Omega_{B}(Y)\ar[r]^{\varphi}&  (A/B)\otimes_{B}P \ar[r]^{\psi}\ar[d]^{}& (A/B)\otimes_{B}Y\ar[r]^{}\ar[d]^{}& 0\\
 &&0&0&
}$$where the third column in the diagram is given by Lemma \ref{3.2}. The exactness of the second row follows from $\Tor_1^B(A,Y)=\Tor_1^B(A,\Omega_B^{n+m+1}(X))\simeq \Tor_{n+m+2}^B(A,X)=0$ since $n$ is a bound on the vanishing of $\Tor^B(A_B,-)$. Thus there is a projective $A$-module $Q$ such that $A\otimes_{B}\Omega_{B}(Y)\simeq \Omega_{A}(A\otimes_{B}Y)\oplus Q$. So we may rewrite the above diagram as follows:
$$\xymatrix{
 &0\ar[d]^{}&0\ar[d]^{}&0\ar[d]^{}& \\
0\ar[r]^{}& \Omega_{B}(Y)\ar[r]^{}\ar[d]^{}&  P\ar[r]^{\pi}\ar[d]^{}& Y\ar[r]^{}\ar[d]^{}& 0\\
0\ar[r]^{}& \Omega_{A}(A\otimes_{B}Y)\oplus Q\ar[r]^{}\ar[d]^{f}&  A\otimes_{B}P \ar[r]^{}\ar[d]^{}& A\otimes_{B}Y\ar[r]^{}\ar[d]^{}& 0\\
0\ar[r]^{}& \Img(\varphi)\ar[r]^{}\ar[d]^{}&  (A/B)\otimes_{B}P \ar[r]^{}\ar[d]^{}& (A/B)\otimes_{B}Y\ar[r]^{}\ar[d]^{}& 0\\
 &0&0&0&
}$$
where $f$ is surjective. Thus we have an exact sequence of $B$-modules
$$ (*)\quad 0\lra \Omega_B(Y)\lra \Omega_A(A\otimes_BY)\oplus Q\lra \Img(\varphi)\lra 0.$$
By Lemma \ref{3.1}, we have the following isomorphism of $A$-modules for any $j\ge 1$:
$$ \Omega_A(A\otimes_B\Omega^{n+j}_B(X)) \simeq \Omega^j_A(A \otimes_B \Omega^{n+1}_B(X)).$$
So, if $j>\pd(_BX)$, then $\Omega^{n+j}_B(X)=0$, and therefore $\Omega_A^j(A \otimes_B \Omega^{n+1}_B(X))=0$ and $\pd(_AA \otimes_B \Omega^{n+1}_B(X))\le m :=\fd(A)$. In particular, it follows again from Lemma \ref{3.1} that $$\Omega_A(A \otimes_B \Omega^{n+m+1}_B(X))\simeq\Omega_A^{m+1}(A\otimes_B\Omega_B^{n+1}(X))=0.$$ Hence the sequence $(*)$ becomes the following  exact sequence
$$ 0\lra \Omega_B^{n+m+2}(X)\lra Q\lra \Img(\varphi)\lra 0.$$
If we take $(s-1)$-th syzygy of this sequence, then we get a new exact sequence of $B$-modules
$$ 0\lra \Omega_B^{n+m+2+s-1}(X)\lra \Omega_B^{s-1}(Q)\oplus P'\lra \Omega_B^{s-1}(\Img(\varphi))\lra 0  $$
with $P'$ a projective $B$-module.
Now, applying the second shifting sequence in Lemma \ref{facts}(2) to this sequence, we obtain the following exact sequence
$$\xymatrix@C=0.5cm{ 0 \ar[r] & \Omega_B^{s}(Q) \ar[r] &\Omega_B^s(\Img(\varphi))\oplus W\ar[r]^{} &\Omega_B^{n+m+s+1}(X) \ar[r] & 0, }$$
where $W$ is a projective $B$-module.

Next, we shall prove that $\Img(\varphi)$ is a $B/I$-module, that is, $I (\Img(\varphi))=0$. In fact, $\varphi$ is induced by inclusion from $\Omega_{B}(Y)$ to $P$, and therefore $\Omega_{B}(Y)\subseteq \rad(_BP)=\rad(B)P$. So an element in $\Img(\varphi)$ is a finite sum of elements of the form $(a+B)\otimes_{B}rp$, where $a\in A$, $r\in \rad(B)$ and $p\in P$. Hence it is enough to show $I \big((a+B)\otimes_{B}rp\big)=0$ in $(A/B)\otimes_BP$. If $r'\in I$, then it follows from $I \, A\, \rad(B)\subseteq B$ that
$$ r'\big((a+B)\otimes_B rp\big)=r'\big((ar+B)\otimes_Bp\big) =(r'a r+B)\otimes_Bp=0\otimes_Bp=0.$$ This shows that $\Img(\varphi)$ is a module over $B/I$. By assumption, the full subcategory $(B/I)\modcat$ of $B\modcat$ is $s$-syzygy-finite. So there is an additive generator $_BN$ for $\Omega_B^s\big((B/I)\modcat\big)$ such that $\Omega_B^s(\Img(\varphi))\in \add(N)$.

Finally, it follows from Lemma \ref{it} that
\begin{align*}
  \pd(_B\Omega_B^{n+m+s+1}(X))&=\Psi(_B\Omega_B^{n+m+s+1}(X))
  \leqslant 1+\Psi\big(\Omega_B^s(Q)\oplus \Omega_{B}^s(\Img(\varphi))\oplus W\big)
   \leqslant 1+\Psi\big(\Omega_B^s(_BA)\oplus N\big).
 \end{align*}
Clearly,  $\Psi\big(\Omega_B^s(_BA)\oplus N \big)$ does not depend upon $X$. As a result, we have
$$\pd(_BX)\leqslant m+n+s+1 +\pd(_B\Omega_{B}^{n+m+s+1}(X))
\leqslant m+n+s+2+\Psi\big(\Omega_B^s(_BA)\oplus N\big)<\infty.$$
That is, $\fd(B)<\infty$. This completes the proof of Theorem \ref{1.1b}. $\square$

\medskip
{\bf Proof of Corollary \ref{1.1}}:

Let $B\subseteq A$ be a right-finite extension of Artin algebras. Suppose that there is an integer $s\ge 0$ such that $\rad^{s}(B)A\,\rad(B)\subseteq B$ and that $B/\rad^{s}(B)$ is representation-finite. Then we define $I:=\rad^s(B)$. Evidently, $I$ satisfies all conditions in Theorem \ref{1.1b}. Hence Corollary \ref{1.1} follows directly from Theorem \ref{1.1b}. $\square$

\smallskip
As an immediate consequence of Corollary \ref{1.1}, we have the following result.

\begin{Koro} Let $B\subseteq A$ be a right-finite extension of Artin algebras. Suppose that there is an integer $s\geq 1$ such that $A\, \rad(B)\subseteq \rad(B)A$, $\rad^{s}(B)$ is a right ideal of $A$, and
$B/\rad^{s-1}(B)$ is  representation-finite. If $\fd(A)<\infty$, then $\fd(B)<\infty$.
\label{Koro1}
\end{Koro}

Recall that the representation dimension of an algebra $A$, introduced by Auslander in \cite{AusRepdim}, is defined by
$$ \rd(A) :=\inf\{\gd(\End_A(A\oplus D(A)\oplus M))\mid M\in A\modcat\}. $$

If we strengthen extensions $B\subseteq A$ as certain special left-finite extensions, we may relax the assumption on $B/I$ in Theorem \ref{1.1b} and get the following result.

\begin{Koro} Let $B\subseteq A$ be a right-bounded extension of Artin algebras with $\pd(_BA)\le 1$, and let $I$ be an ideal of $B$ such that $I\, A\, \rad(B)\subseteq B$ and $\rd(B/I)\le 3$. If $\fd(A)<\infty$, then $\fd(B)<\infty$.
\label{Koro1a}
\end{Koro}

{\it Proof.} As in the proof of Theorem \ref{1.1b}, we get an exact sequence of $B$-modules
$$\xymatrix@C=0.5cm{ 0 \ar[r] & \Omega_B(Q) \ar[r] &\Omega_B(\Img(\varphi))\oplus W\ar[r]^{} &\Omega_B^{n+m+2}(X)\oplus V \ar[r] & 0, }$$
where $Q$ is a projective $A$-module, $W$ and $V$ are projective $B$-modules. Since $\rd(B/I)\le 3$ and $\Img(\varphi)$ is a $B/I$-module, we can find a $B/I$-module $U$ such that $\gd(\End_{B/I}(U))=\rd(B/I)$, and an exact sequence $$0\lra U_1\lra U_0\lra \Img(\varphi)\lra 0$$
such that $U_i\in \add(U)$ by Lemma \ref{repdim}. Thus we can form the following exact commutative diagram of $B$-modules:
$$\xymatrix{
 &0\ar[d]^{} &0\ar[d]^{}& & \\
 & \Omega_B(U_1)\ar@{=}[r]\ar[d]^{}&  \Omega_B(U_1) \ar[d]^{}& & \\
 0\ar[r]^{}& K\ar[r]^{}\ar[d]^{}& \Omega_B(U_0)\oplus W' \ar[r]^{}\ar[d]^{}& \Omega_B^{n+m+1}(X)\oplus V\ar[r]^{}\ar@{=}[d]& 0\\
0\ar[r]^{}& \Omega_{B}(Q)\ar[r]^{}\ar[d]^{}&  \Omega_B(\Img(\varphi))\oplus W\ar[r]^{}\ar[d]^{}& \Omega_B^{n+m+1}(X)\oplus V\ar[r]^{}& 0\\
 &0&0&&
}$$
where $W'$ is a projective $B$-module. It follows from $\pd(_BA)\le 1$ that $\ext^1_B(\Omega_B(Q),\Omega_B(U_1))=\ext^2_B(Q,\Omega_B(U_1))=0$, and therefore the first column of the above diagram splits. Thus $K\simeq \Omega_B(Q\oplus U_1)$, and consequently, we have the following exact sequence
$$\xymatrix@C=0.5cm{ 0 \ar[r] & \Omega_B(Q\oplus U_1) \ar[r] &\Omega_B(U_0)\oplus W'\ar[r]^{} &\Omega_B^{n+m+2}(X)\oplus V \ar[r] & 0. }$$
Now, the same argument of the proof of Theorem \ref{1.1b} will lead to Corollary \ref{Koro1a}. $\square$

\subsection{Proof of Theorem \ref{1.3}}

Before starting with the proof of Theorem \ref{1.3}, we first show the following useful lemma which supplies a non-trivial way to lift $B$-modules
to $A$-modules. This observation extends \cite[Lemma 0.1]{I} and is probably not previously known in the literature. Also, note that the lemma holds true for arbitrary rings and  torsionless modules.

\begin{Lem}
Let $B \subseteq A$ be an extension of Artin algebras and $I$ be an ideal of $B$ such that $I$ is also a left ideal of $A$. Then, for any torsionless $B$-module $X$, its submodule $IX$ admits an $A$-module structure and is actually a torsinless $A$-module.
\label{Lem3}
\end{Lem}

{\it Proof.} Let $X\hookrightarrow P$ be an inclusion with $P$ a projective $B$-module. Then we get an exact sequence $$0\lra X\inja P\lra Y\lra 0$$ with $Y$ the cokernel of the inclusion, and can form the following exact commutative diagram:

$$\xymatrix{
0\ar[r] &\Tor_1^B(I,Y)\ar[r]^{\psi}&I\otimes_BX\ar[r]^{\eta}&I\otimes_B P\ar[d]^{\alpha}\ar[r]&I\otimes_B Y\ar[d]^{\beta}\ar[r]& 0\\
 &  &    &IP\ar[r]&IY\ar[r] & 0,}$$
where $\alpha$ and $\beta$ are the multiplication maps. Clearly, $\alpha$ is an isomorphism since $P$ is a projective $B$-module. It is easy to see that $\Img(\eta)\simeq \Img(\eta\alpha)= IX$. Thus we have an exact sequence of $B$-modules:
\[
\xymatrix{
0\ar[r]&\Tor_1^B(I,Y)\ar[r]^{\psi}&I\otimes_BX\ar[r]^{\eta\alpha}&IX\ar[r]&0.
}
\]
Since $I$ is a left ideal of $A$, we know that $\Tor_1^B(I,Y)$ and $I\otimes_BX$ are $A$-modules and the injective homomorphism $\psi$ is, in fact,
a homomorphism of $A$-modules. So the $B$-module $IX$, as the cokernel of $\psi$, is endowed with an $A$-module structure which is induced from the
$A$-module $I\otimes_BX$. Evidently, the $A$-module $IX$ can be regarded as an $A$-submodule of the projective $A$-module $_AA\otimes_BP$. $\square$

\medskip
Remark that Lemma \ref{Lem3} may be false if $X$ is not assumed to be torsionless. For a counterexample, we refer the reader to \cite[Erratum]{I}.

\medskip
{\bf Proof of Theorem \ref{1.3}}:

Let $X$ be a $B$-module with $\pd(_BX)< \infty$. Then we take a projective cover $P_1\ra \Omega_B(X)$ of $\Omega_B(X)$ and get two canonical exact sequences of torsionless $B$-modules:
\[0 \lra \Omega^2_B(X) \lra P_1 \lra \Omega_B(X) \lra 0,\]
\[0 \lra \Omega^2_B(X) \lra \rad_B( P_1) \lra \rad_B(\Omega_B(X)) \lra 0,\]
which can be used to construct the following exact commutative diagram of $B$-modules:
\[
\xymatrix{
&&0\ar[d] &0\ar[d] &0\ar[d] &\\
&0\ar[r]&Z\ar[d]^{\iota_1}\ar[r]&I \, \rad_B(P_1)\ar[d]^{\iota_2}\ar[r]^{\varphi}&I \, \rad_B(\Omega_B(X))\ar[d]^{\iota_3}\ar[r]&0\\
&0\ar[r]&\Omega^2_B(X)\ar[d]\ar[r]&\rad_B( P_1)\ar[d]\ar[r]&\rad_B(\Omega_B(X))\ar[d]\ar[r]&0\\
&0\ar[r]&T\ar[d]\ar[r]&\rad_B( P_1)/I \, \rad_B(P_1)\ar[d]\ar[r]&\rad_B(\Omega_B(X))/I \, \rad_B(\Omega_B(X))\ar[d]\ar[r]&0\\
&&0&0&0&
}
\]
where $Z$ is the kernel of $\varphi$ and where $T$ is the cokernel of $\iota_1$. Thus we have an exact sequence of $B$-modules
$$(**)\quad 0\lra Z\lra \Omega_B^2(X)\lra T\lra 0 $$
with $T$ a $B/I$-module. Note that $Z$ has an $A$-module structure: In fact, $\varphi$ is the composite of $I \, \rad(B)\otimes_BP_1\ra I \, \rad(B)\otimes_B\Omega_B(X)$
with the multiplication map $I \, \rad(B)\otimes_B\Omega_B(X)\ra I \, \rad(B)\Omega_B(X)=I \, \rad_B(\Omega_B(X))$. Since $I \, \rad(B)$ is a left ideal of $A$ and $\Omega_B(X)$ is a
torsionless $B$-module, we see that $I \, \rad_B(\Omega_B(X))$ is an $A$-module by Lemma \ref{Lem3}. It turns out that $\varphi$ becomes a homomorphism of $A$-modules, and therefore
its kernel $Z$ is an $A$-module.

Since $P_1$ is a projective $B$-module, we have the following inclusions of $A$-modules:
$$Z \hookrightarrow I \, \rad_B(P_1) \simeq I \, \rad(B) \otimes_B P_1 \hookrightarrow A \otimes_B P_1.$$
Let $W$ be the cokernel of the inclusion $Z \hookrightarrow A\otimes_B P_1$. Then
there is a projective $A$-module $Q$ such that $Z \simeq \Omega_A(W)\oplus Q$. Thus the exact sequence ($**$) can be rewritten as
$$ 0 \ra \Omega_A(W)\oplus Q \ra \Omega^2_B(X) \ra T \ra 0.$$
Suppose that the full subcategory $(B/I)\modcat$ of $B\modcat$ is $s$-syzygy-finite for some $s\ge 1$. Then we take the $(s-1)$-th syzygy of the above sequence and get the following exact sequence of $B$-modules:
$$ 0 \ra \Omega_B^{s-1}(\Omega_A(W)\oplus Q) \ra \Omega^{s+1}_B(X)\oplus P' \ra \Omega_B^{s-1}(T) \ra 0.$$
with $P'$ a projective $B$-module.
By a syzygy shifting, we finally get an exact sequence of the following form:
$$\xymatrix@C=0.5cm{
0 \ar[r] & \Omega_{B}^s(T) \ar[r]^{} & \Omega_B^{s-1}(\Omega_A(W)\oplus Q)\oplus P \ar[r]^{} & \Omega^{s+1}_B(X) \ar[r] & 0, }$$where $P$ is a projective $B$-module. Since $A$ is torsionless-finite and the full subcategory $(B/I)\modcat$ of $B\modcat$ is $s$-syzygy-finite, there is an $A$-module $M$ and a $B$-module $N$ such that $\Omega_A(W)\in \add(_AM)$ and $\Omega_B^s(T)\in \add(_BN)$. Here $M$ and $N$ do not depend upon $X$. Thus the Igusa-Todorov function yields the following estimation:
  \begin{align*}
  \pd(_BX) &\leq s+1+\pd(_B\Omega^{s+1}_B(X))\\
  &\leq s+1+1+\Psi\big(\Omega_B^s(T)\oplus \Omega_B^{s-1}(\Omega_A(W) \oplus Q )\oplus P\big) \\
  &\leq s+2+\Psi(N\oplus \Omega_B^{s-1}(M \oplus {}_BA)).
  \end{align*}
Note that $s$ and $\Psi(N\oplus \Omega_B^{s-1}(M \oplus {}_BA))$ are independent of the $B$-module $X$. Hence $\fd(B)<\infty$. This finishes the proof of Theorem \ref{1.3}. $\square$

\medskip
Remarks. (1) The above proof shows that, in Theorem \ref{1.3}, we may replace ``$A$ is torsionless-finite"  by ``the subcategory $\Omega_B^{s-1}\big(\Omega_A(A\modcat)\big)$ is of finite type".

 (2) Given an Artin $R$-algebra $B$ and an ideal $I$ in $B$, if $B$ is a free $R$-module (for instance, $R$ is a field), then there is a recipe for getting an extension $B\subseteq A$ of Artin $R$-algebras such that $I$ is a left ideal in $A$. In fact, since $B$ is a free $R$-module of finite rank, we can embed $B$ into a full $n\times n$ matrix algebra $\Lambda: =M_n(R)$ over $R$ and define $A:=\{a\in \Lambda \mid aI\subseteq I\}$. Evidently, the Artin $R$-algebra $A$ contains $B$ and makes $I$ into a left ideal.

\medskip
{\bf Proof of Corollary \ref{1.2}}:

Suppose that $B\subseteq A$ is an extension of Artin algebras such that $\rad^{s}(B)$ is a left ideal of $A$ and that $B/\rad^{s-1}(B)$ is representation-finite for some integer $s\geq 1$. If $I:=\rad^{s-1}(B)$, then $I$ fulfills the conditions in Theorem \ref{1.3}. So Corollary \ref{1.2} follows immediately from Theorem \ref{1.3}. $\square$

\medskip
The above proof of Theorem \ref{1.3} also implies the following conclusion.

\begin{Koro}
Let $B \subseteq A$ be a left-finite extension of Artin algebras such that $\rad^2(B)$ is a left ideal of $A$.

 $(1)$ If $A$ is torsionless-finite, then $\fd(B) < \infty$.

 $(2)$ If $\gd(A)<\infty$, then $\fd(B) < \infty$.
\label{Koro4}
\end{Koro}

{\it Proof.} (1) is clear from Corollary \ref{1.2}.

(2) We keep the notation in the above proof of Theorem \ref{1.3}. By Lemma \ref{facts}(1), $\pd(_BZ)\le \pd(_BA)+\gd(A)$. This shows $\pd(_BT)<\infty$ for the semisimple module $T$. Hence, the above argument yields $\fd(B)<\infty$. $\square$

\medskip
Now, we consider extensions $B\subseteq A$ with $\rad^l(B)=\rad^l(A)$. For $l=1$, it is known from \cite[Theorem 1.1(1)]{XiXu} that $\fd(B)<\infty$ if $\fd(A)<\infty$. For $l\ge 2$, we have the following variation of Corollary \ref{1.2}.

\begin{Koro}
Let $B \subseteq A$ be an extension of Artin algebras such that $\rad^l(B)=\rad^l(A)$ for an integer $l\geq 2$ and that both $A/\rad^{l-1}(A)$ and $B/\rad^{l-1}(B)$ are representation-finite. If $\gd(A) \le 2$, then $\fd(B) < \infty$.
\label{Koro7}
\end{Koro}

{\it Proof.} Let $X$ be a $B$-module with $\pd(_{B}X)< \infty$, and let $\cdots \ra P_2\ra P_1\ra P_0\ra X\ra 0$ be a minimal projective resolution of $_BX$. Then, as in the proof of Theorem \ref{1.3}, we can construct an exact sequence
\[
(**)\quad 0 \lra Z \lra \Omega^2_B(X) \lra T \lra 0,
\]
 of $B$-modules, which induces another exact sequence of $A$-modules $$0 \lra Z \lra \rad^l_B(P_1) \lra \rad^l_B(\Omega_B(X)) \lra 0,$$
where these $A$-module structures are due to Lemma \ref{Lem3} and where $T$ is a module over $B/\rad^{l-1}(B)$. Furthermore, we also have the following exact sequence of $A$-modules:$$0 \lra \Omega^2_A(A\otimes_BX)\oplus P \lra A\otimes_B P_1 \lra A\otimes_B P_0 \lra A\otimes_BX \lra 0$$with $P$ a projective $A$-module. Now, we construct the following exact commutative diagram of $A$-modules:
\[
\xymatrix{
&0\ar[d]&0\ar[d]&0\ar[d]&&\\
0\ar[r]&Z\ar[d]\ar[r]&\rad^l_B(P_1)\ar[d]\ar[r]&\rad^l_B(\Omega_B(X))\ar[d]\ar[r]&0&\\
0\ar[r]&\Omega^2_A(A\otimes_BX)\oplus P\ar[d]\ar[r]&A\otimes_B P_1 \ar[d]\ar[r]&A\otimes_B P_0\ar[d]\ar[r]&A\otimes_BX\ar@{=}[d]\ar[r]&0\\
0\ar[r]&S_1\ar[r]\ar[d]&(A/\rad^l(B))\otimes_BP_1\ar[r]\ar[d]&S_3\ar[r]\ar[d]&A\otimes_BX\ar[r]&0\\
&0&0&0&&
}
\]
Since $\rad^l(B)=\rad^l(A)$, we see that $S_1$ is a module over $A/\rad^l(A)$ and that $(A/\rad^l(B))\otimes_BP_1$, as an $A/\rad^{l}(A)$-module, is projective. Thus $S_1$ is a torsionless module over $A/\rad^l(A)$.
By assumption, $A/\rad^{l-1}(A)$ is representation-finite, this implies that $A/\rad^l(A)$ is torsionless-finite (see \cite{AusRepdim, ringel}). So there is a module $M$ over $A/\rad^{l}(A)$ such that $S_1\in \add(M)=\Omega^1_{A/\rad^l(A)}\big((A/\rad^l(A))\modcat\big)$. Since $\gd(A) \le 2$, we see that the $A$-module $\Omega^2_A(A\otimes_BX)$ is projective and $Z\simeq \Omega_A(S_1)\oplus Q$ with $Q$ a projective $A$-module. Consequently, $Z \in \add(\Omega_A(M)\oplus {}_AA)$.

Since $B/\rad^{l-1}(B)$ is representation-finite, we may find an additive generator $N$ for $\big(B/\rad^{l-1}(B)\big)\modcat$ such that $\Omega_B(T)\in \add(\Omega_B(N))$.
Now, applying Lemma \ref{facts} to $(**)$, we have an exact sequence:
\[
0 \lra \Omega_B(T) \lra Z\oplus P' \lra \Omega^2_B(X) \lra 0
\]
with $P'$ a projective $B$-module. By the Igusa-Todorov function (see Lemma \ref{it}), we have the following inequalities:
\begin{align*}
\pd({}_{B}X) &\le 2+\pd({}_{B}\Omega^2_B(X))
\le 2+\Psi(Z\oplus P'\oplus \Omega_B(T))
\le 2+\Psi\big(_B\Omega_A(M)\oplus {}_BA \oplus \Omega_B(N)\big).
\end{align*}
So $\fd(B) < \infty$. $\square$

\bigskip
In Corollary \ref{Koro7}, we take $l=2$ and get the following corollary.

\begin{Koro}
Let $B \subseteq A$ be an extension of Artin algebras with $\rad^2(B)=\rad^2(A)$. If $\gd(A) \le 2$, then $\fd(B) < \infty$.
\label{Koro5}
\end{Koro}

Corollary \ref{Koro5} can be used to re-obtain a main result in \cite{Green} on algebras with vanishing radical cube.

\begin{Koro}\emph{\cite{Green}} Suppose that $B$ is a finite-dimensional $k$-algebra of the form $B=kQ/I$ with vanishing radial cube,
where $Q$ is a quiver and $I$ is an admissible ideal in the path algebra $kQ$ of $Q$ over the field $k$. Then $\fd(B)<\infty$.
\label{koro6}
\end{Koro}

{\it Proof.} Since $\rad^3(B)=0$, we may write $B$ as $B=kQ_0\oplus kQ_1\oplus kQ_2$, where $Q_0$ and $Q_1$ are the sets of vertices and arrows of $Q$, respectively, and $Q_2$ is a set of $k$-linearly independent paths of length two in $B$. Then $B$ is embedded canonically into a triangular matrix algebra $A:=\begin{pmatrix} kQ_0        &  &  \\
  kQ_1     & kQ_0   &         \\  kQ_2 & kQ_1 & kQ_0 \\
\end{pmatrix}$ by sending $b=b_0+b_1+b_2$ to $\begin{pmatrix} b_0        &  &  \\
  b_1     & b_0   &         \\  b_2 & b_1 & b_0 \\
\end{pmatrix}$, where $b_i\in kQ_i$ for $ 0\le i \le 2$.
We can check that $\gd(A)\le 2$ and $\rad^2(B)=\rad^2(A)$. Thus Corollary \ref{koro6} follows from Corollary \ref{Koro5}. $\square$

\medskip
For a chain of extensions, we have the following result which extends \cite[Theorem 4.5]{I} slightly.
\begin{Koro}
Let $C \subseteq B \subseteq A$ be a chain of extensions of Artin algebras such that $\rad(C)$ is a left ideal of $B$ and $\rad^l(B)$ is a left ideal of $A$ for some integer $l\geq 1$. If $A$ is torsionless-finite and $B/\rad^l(B)$ is representation-finite, then $\fd(C)< \infty$.
\label{Koro6}
\end{Koro}

{\it Proof.} Let $X$ be a $C$-module with $\pd(_{C}X)< \infty$. By Lemma \ref{Lem4}, there is a $B$-module $Y$ and a projective $B$-module $P$ such that $\Omega^2_C(X) \simeq \Omega_B(Y)\oplus P$ as $B$-modules. Clearly, the $C$-modules $P$ and $\Omega_B(Y)$ have finite projective dimensions. It follows from Lemma \ref{facts}(2) and the exact sequence
\[
0 \lra \rad^l_B(\Omega_B(Y)) \lra \Omega_B(Y) \lra \Omega_B(Y)/\rad^l_B(\Omega_B(Y)) \lra 0
\]
that the following sequence of $C$-modules
\[
0 \lra \Omega_C\big(\Omega_B(Y)/\rad^l_B(\Omega_B(Y))\big) \lra \rad^l_B(\Omega_B(Y))\oplus Q \lra \Omega_B(Y) \lra  0
\]
is exact, where $Q$ is a projective $C$-module.
By Lemma \ref{Lem3}, $\rad^l_B(\Omega_B(Y))$ is a torsionless $A$-module. Clearly, $\Omega_B(Y)/\rad^l_B(\Omega_B(Y))$ is a module over $B/\rad^l(B)$.
Since $A$ is torsionless-finite and $B/\rad^{l}(B)$ is representation-finite, we may assume that $M$ and $N$ are additive generators for $\Omega(A\modcat)$ and $(B/\rad^{l}(B))\modcat$, respectively. Thus
$\rad^l_B(\Omega_B(Y))\in \add(M)$ and $\Omega_B(Y)/\rad^l_B(\Omega_B(Y))\in \add(_BN)$.

Now, we use Lemma \ref{it} to give an upper bound for the projective dimension of $_CX$:
\begin{align*}
\pd({}_{C}X) &\le 2+\pd({}_{C}\Omega^2_C(X))\\
&= 2+\pd({}_{C}\Omega_B(Y)\oplus {}_CP)\\
&\le 2+\Psi\bigg(\Omega_C\big(\Omega_B(Y)/\rad^l_B(\Omega_B(Y))\big)\oplus \rad^l_B(\Omega_B(Y))\oplus Q\oplus P\bigg)\\
&\le 2+\Psi(\Omega_C(N)\oplus M \oplus {}_CB).
\end{align*}
So we have $\fd(C) < \infty$. $\square$

\medskip
Next, we point out that there are lots of right-finite extensions $B\subset A$ satisfying the condition $\rad(B)A\subset B$, and therefore $\rad(B) A\, \rad(B)\subseteq B$. This is the case of Corollary \ref{1.1} for $s=1$. Let us exhibit one such example.

Suppose that $\Lambda$ is an Artin algebra. We define $A :=M_2(\Lambda)$, the algebra of $2\times 2$ matrices over $\Lambda$, and
$$B:=\begin{pmatrix}
  \Lambda        &  \rad(\Lambda)  \\
  \Lambda        &  \Lambda        \\
 \end{pmatrix}.$$
Then $A_B$ and $_BA$ are projective and the extension $B\subseteq A$ is both right- and left-finite. An easy calculation shows that
$$\rad(B)=\begin{pmatrix}
  \rad(\Lambda)        &  \rad(\Lambda)        \\
  \Lambda              &  \rad(\Lambda)
\end{pmatrix},\quad
\rad(A)=\begin{pmatrix}
  \rad(\Lambda)        &  \rad(\Lambda)    \\
  \rad(\Lambda)        &  \rad(\Lambda)
\end{pmatrix}\quad \mbox{and} \quad \rad(B)A=\begin{pmatrix}
  \rad(\Lambda)        &  \rad(\Lambda)        \\
  \Lambda              &  \Lambda
\end{pmatrix},$$
and that $\rad(B)$ is neither a left nor a right ideal in $A$. But we have $\rad(A)\subsetneq\rad(B)\subseteq\rad(B)A\subseteq B$, as desired. If $\fd(\Lambda)<\infty$, then $\fd(B)<\infty$ by Corollary \ref{1.1} for $s=1$.
In this example, if $\rad(\Lambda)\ne 0$, then $\rad^2(B)$ is neither a left nor a right ideal in $A$, but $\rad^3(B)=\rad^2(A)$ is an ideal in $A$. Moreover, if $\fd(\Lambda)<\infty$, then we can get $\fd(B)<\infty$ alternatively by Corollary \ref{Koro1a} since $\rad^2(B) A\,\rad(B)\subseteq B$ and since $B/\rad^2(B)$ always has representation dimension at most $3$ by a result of Auslander (see \cite[Chapt. III, Sec. 5]{AusRepdim}).

\medskip
Finally, we display an example to show how our results developed in this paper can be applied to decide whether certain algebras have finite finitistic dimension.
The example shows also that the method of controlling finitistic dimensions by extension algebras seems to be useful.

Let $A$ be an algebra (over a field) given by the following quiver
\[
\xymatrix{
5\bullet     & 2\bullet\ar[l]_{\lambda}& 3\bullet\ar[l]_{\epsilon}      &1\bullet\ar[l]_{\xi}&4\bullet\ar[l]_{\beta} &\bullet 6\ar[l]_{\alpha}\\
}
\]
with one relation: $\alpha\beta\xi\epsilon\lambda=0$. Clearly, this algebra is representation-finite.
Now, let $B$ be the subalgebra of $A$ generated by $\{e_1, \; e_{2'}:=e_2+e_4+e_5, \; e_{3'}:=e_3+e_6, \lambda, \beta, \alpha+\epsilon, \gamma:=\xi\epsilon, \delta:=\beta\xi \}$, where $e_i$ is the primitive idempotent element of $A$ corresponding to the vertex $i$. Then $B$ is given by the following quiver
\[
\xymatrix{
&1\bullet\ar@<1ex>[r]^{\gamma}& \ar@(ur,ul)[]|{\lambda}\;\;\bullet{2'}\;\ar@<1ex>[l]^{\beta}\;\;\ar@<1ex>[r]^{\delta}&\bullet{3'}\ar@<1ex>[l]^{\quad\alpha+\epsilon}&
}
\]
with relations: $\beta\gamma=\delta(\alpha+\epsilon)$, $\gamma\beta=\gamma\delta=\lambda^2=\lambda\beta=\lambda\delta=(\alpha+\epsilon)\beta\gamma\lambda=0$, and the Loewy structures of the indecomposable projective $B$-modules are as follows:

\vspace{0.5cm}
{\unitlength=0.5cm
\special{em:linewidth 0.4pt}
\linethickness{0.4pt}
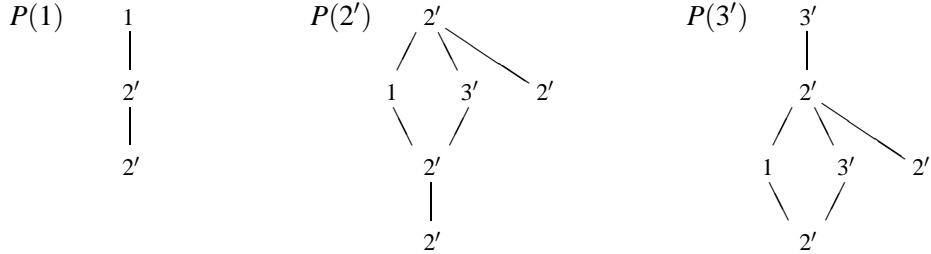
\begin{center}
\begin{picture}(20,5.2)
\put(-2,2){{\footnotesize $2'$}}
\put(-2,4){{\footnotesize $2'$}}
\put(-2,6){{\footnotesize $1$}}
               \put(6,6){{\footnotesize $2'$}}
\put(5,4){{\footnotesize $1$}}   \put(7,4){{\footnotesize $3'$}}  \put(9,4){{\footnotesize $2'$}}
                    \put(6,2){{\footnotesize $2'$}}
                    \put(6,0){{\footnotesize $2'$}}
         \put(16,6){{\footnotesize $3'$}}
         \put(16,4){{\footnotesize $2'$}}
\put(15,2){{\footnotesize $1$}} \put(17,2){{\footnotesize $3'$}} \put(19,2){{\footnotesize $2'$}}
              \put(16,0){{\footnotesize $2'$}}
\put(-1.8,5.8){\line(0,-1){1}}\put(-1.8,3.8){\line(0,-1){1}}

\put(5.8,5.8){\line(-1,-2){0.5}}\put(6.4,5.8){\line(1,-2){0.5}}
\put(6.6,5.9){\line(3,-2){2.2}}
\put(5.2,3.8){\line(1,-2){0.5}}

\put(7.2,3.8){\line(-1,-2){0.5}}\put(6.2,1.8){\line(0,-1){1}}

\put(16.2,5.8){\line(0,-1){1}}
\put(15.8,3.8){\line(-1,-2){0.5}}\put(16.4,3.8){\line(1,-2){0.5}}
\put(16.6,3.9){\line(3,-2){2.2}}
\put(15.2,1.8){\line(1,-2){0.5}}\put(17.2,1.8){\line(-1,-2){0.5}}

\put(-5,6){$P(1)$}\put(3,6){$P(2')$}\put(13,6){$P(3')$}
\end{picture}
\end{center}}\vspace{-0.3cm}\noindent It is not difficult to see that $B$ is representation-infinite and all simple $B$-modules have infinite projective dimension.
Thus the length $\ell\ell^{\infty}(B)$, defined in \cite{HLM}, is just the Loewy length of $B$. Moreover, the algebra $B$ is neither monomial nor
radical-cube vanishing nor standardly stratified nor special biserial. Note that $B/\rad^3(B)$ is representation-infinite and that all of $\pd(_BA)$,
$\pd(A_B)$, $\pd(_B\rad^i(B))$
and $\pd(\rad^i(B)_B)$ for $1\le i\le 3$ are infinite. So it is not clear that $\rd(B)\le 3$. Though $B$ is embedded into $A$ of representation dimension $2$, the result \cite[Theorem 4.2]{II} cannot be applied because $\rad(B)$ is neither a left nor a right ideal in $A$. But we can verify that $\rad^2(B)$ is an ideal of $A$, and therefore $\fd(B)< \infty$ by Corollary \ref{1.2}.

\medskip
{\bf Acknowledgement.} The paper was partially revised during a visit of the corresponding author Changchang Xi to the University of Stuttgart, Germany, from June to August in 2015.
He enjoyed the stay in Stuttgart very much and would like to thank Steffen Koenig for his invitation and hospitality. The research work is partially supported by BNSF and NNSF (KZ201410028033, 11331006).

{\footnotesize
\bigskip

Chengxi Wang

School of Mathematical Sciences, Beijing Normal University,
100875 Beijing, People's Republic of China

{\tt Email: chxwang66@mail.bnu.edu.cn}

\bigskip
Changchang Xi

School of Mathematical Sciences, BCMIIS, Capital Normal University, 100048
Beijing, People's Republic of  China

{\tt Email: xicc@cnu.edu.cn}}

\medskip
{\footnotesize First version: May 16, 2015, Revised: July 27, 2015}

\end{document}